\documentclass[11pt,reqno]{amsart}

\usepackage{amssymb, amsmath, amsfonts, latexsym}
\usepackage{enumerate}
\usepackage[notref,notcite]{showkeys}

\newtheorem{thm}{Theorem}[section]
\newtheorem{cor}[thm]{Corollary}
\newtheorem{lem}[thm]{Lemma}
\newtheorem{prop}[thm]{Proposition}
\newtheorem{example}[thm]{Example}
\newtheorem{remarks}[thm]{Remarks}
\newtheorem{defn}[thm]{Definition}

\numberwithin{equation}{section}

\def\ff{\frac}
\def\ss{\sqrt}
\def\BB{\mathbb B}

\def\DD{\Delta}

\def\<{\langle}
\def\>{\rangle}

\def\nnn{\nabla}
\def\d{\text{\rm{d}}}
\def\aa{\alpha}

\def\E{\mathcal E}

\def\beg{\begin}
\def\beq{\begin{equation}}

\def\Ric{\text{\rm{Ric}}}
\def\Hess{\text{\rm{Hess}}}

\def\e{\text{\rm{e}}}

\def\Ric{\text{\rm{Ric}}}

\title[Transport/Fisher Information inequalities II]{{\bf Transportation-information inequalities for Markov processes (II) : relations
with other functional inequalities}}

\author{Arnaud Guillin}
\address{Arnaud Guillin. Ecole Centrale de Marseille et LATP, Centre de Math\'ematiques et Informatique. Technop\^ole de Ch\^ateau-Gombert. 13453 Marseille, France}
\email{guillin@cmi.univ-mrs.fr}

\author{Christian L\'eonard}
\address{Christian L\'eonard. Modal-X, Universit\'e Paris 10. B\^atiment G, 200 avenue de la R\'epublique.92001 Nanterre, France and CMAP, \'Ecole Polytechnique. 91128 Palaiseau, France}
\email{christian.leonard@polytechnique.fr}

\author{Feng-Yu Wang}
\address{Feng Yu Wang. Department of Mathematics, Swansea University, Singleton Park, SA2 8PP, Swansea, UK}
\email{f.y.wang@swansea.ac.uk}

\author{Liming Wu}
\address{Liming Wu. Laboratoire de Math\'ematiques Appliqu\'ees, CNRS-UMR 6620, Universit\'e Blaise Pascal, 63177 Aubi\`ere, France.
And Department of Mathematics, Wuhan University, 430072 Hubei,
China} \email{Li-Ming.Wu@math.univ-bpclermont.fr}

\date{First Version: March 2008}

\newcommand{\dd}{\mathbb{D}}
\newcommand{\ee}{\mathbb{E}}

\newcommand{\nn}{\mathbb{N}}
\newcommand{\rr}{\mathbb{R}}
\newcommand{\pp}{\mathbb{P}}

\def\AA{\mathcal A}
\def\BB{\mathcal B}
\def\CC{\mathcal C}
\def\DD{\mathcal D}

\def\EE{\mathcal E}

\def\LL{\mathcal L}
\def\MM{\mathcal M}
\def\NN{\mathcal N}

\def\VV{\mathcal V}
\def\XX{\mathcal X}

\def\vep{\varepsilon}
\def\<{\langle}
\def\>{\rangle}

\def\beq{\begin{equation}}
\def\neq{\end{equation}}

\def\bdef{\begin{defn}}
\def\ndef{\end{defn}}

\def\bthm{\begin{thm}}
\def\nthm{\end{thm}}

\def\bprop{\begin{prop}}
\def\nprop{\end{prop}}

\def\brmk{\begin{remarks}}
\def\nrmk{\end{remarks}}

\def\bexa{\begin{example}}
\def\nexa{\end{example}}

\def\blem{\begin{lem}}
\def\nlem{\end{lem}}

\def\bcor{\begin{cor}}
\def\ncor{\end{cor}}

\def\bexe{\begin{exe}}
\def\nexe{\end{exe}}

\def\bprf{\begin{proof}}
\def\nprf{\end{proof}}

\def\dsp{\displaystyle}

\def\bdes{\begin{description}}
\def\ndes{\end{description}}

\def\beg{\begin}

\begin{document}
\maketitle

\begin{abstract}  We continue our investigation on the
transportation-information inequalities $W_pI$ for a symmetric
markov process, introduced and studied in \cite{GLWY}. We prove that
$W_pI$ implies the usual transportation inequalities $W_pH$, then
the corresponding concentration inequalities for the invariant
measure $\mu$.
 We give also a direct proof that the spectral gap in
the space of Lipschitz functions for a diffusion process implies
$W_1I$ (a result due to \cite{GLWY}) and a Cheeger type's
isoperimetric inequality. Finally we exhibit relations between
transportation-information inequalities and a family of functional
inequalities (such as $\Phi$-log Sobolev or $\Phi$-Sobolev).
\end{abstract}

\vskip20pt\noindent {\bf keywords:}  Wasserstein distance; entropy; Fisher information;
transport-information inequality; deviation inequality.


\section{Introduction}

Let $(\XX,d)$ be a complete and separable metric space (say Polish)
and $\mu$ a given probability measure on $(\XX ,\BB)$ where $\BB$ is
the Borel $\sigma$-field. Let $(X_t)_{t\ge0}$ be a $\mu$-symmetric
ergodic conservative  Markov process valued in $\XX$, with
transition semigroup $(P_t)$ (which is symmetric on $L^2(\mu)$), and
Dirichlet form $(\EE(\cdot,\cdot), \dd(\EE))$ where $\dd(\EE)$ is
the domain of $\EE$ in $L^2(\mu):=L^2(\XX ,\BB,\mu)$. Here the
ergodicity means simply : for $g\in\dd(\EE)$, $\EE(g,g)=0$ iff
$g=c$.

For $1\le p<+\infty$ fixed and for any probability measure $\nu$ on
$\XX$ (written as $\nu\in\MM_1(\XX )$), consider

\bdes
\item[(i)] {\bf $L^p$-Wasserstein distance between $\nu$ and $\mu$}:

\beq W_p(\nu,\mu):=\inf_{\pi\in \CC(\nu,\mu)} \iint_{E^2} d^p(x,y)
\pi(dx,dy) \neq where $\CC(\nu,\mu)$ are the set of all couplings of
$(\nu,\mu)$, i.e., probability measures $\pi$ on $E\times E$ such
that  $\pi(A\times E)=\nu(A)$ and $\pi(\XX \times A)=\mu(A)$ for all
$A\in \BB$.

\item[(ii)]  {\bf Relative entropy or Kullback's information of $\nu$ w.r.t. $\mu$}

\beq\label{def-entropy-a} H(\nu|\mu):=\begin{cases} \int
\frac{d\nu}{d\mu} \log \frac{d\nu}{d\mu} d\mu,\ \ &\text{ if }\
\nu\ll
\mu;\\
+\infty, &\text{ otherwise.}
 \end{cases}
\neq
\item[(iii)] {\bf The Fisher information of $\nu$ w.r.t. $\mu$}:

\beq\label{Fisher} I(\nu|\mu) :=
  \begin{cases}
    \EE(\sqrt{f},\sqrt{f} ).\ & \text{ if }\nu=f \mu, \sqrt{f}\in\dd(\EE), \\
    +\infty & \text{otherwise}.
  \end{cases}
\neq  \ndes

The usual transport inequalities $W_pH$, introduced and studied by
K. Marton \cite{Mar96} and M. Talagrand \cite{Tal96a} mean that

\beq W_p(\nu,\mu)^2 \le 2C H(\nu|\mu), \ \ \forall \nu\in M_1(\XX
).\tag{$W_pH(C)$} \neq Its study is very active: see Bobkov-G\"otze
\cite{BG99}, Otto-Villani \cite{OVill00}, Bobkov-Gentil-Ledoux
\cite{BGL01}, Djellout-Guillin-Wu \cite{DGW03} and references
therein. Furthermore Gozlan-L\'eonard \cite{GL} consider the
following generalized transportation cost from $\nu$ to $\mu$:
$$
T_\VV (\nu,\mu):=\sup\{\nu(u)-\mu(v);\ (u,v)\in\VV\}
$$
($\mu(u):=\int_E ud\mu$) where $\VV$ is some given family of
$(u,v)\in (b\BB)^2$ so that

\begin{itemize}
    \item[(A1)] $u\le v$ for all $(u,v)\in \VV$ (or
equivalently $T_\VV(\nu, \nu)\le 0$ for all $\nu\in M_1(\XX)$);
    \item[(A2)] For all $\nu_1,\nu_2\in M_1(\XX)$, there exists $(u,v)\in
\VV$ such that $\int u\, d\nu_1 - \int v\,d\nu_2\ge 0$ (or
equivalently $T_\VV(\nu_1, \nu_2)\ge 0$ for all $\nu_1, \nu_2\in
M_1(\XX)$).
\end{itemize}
And they introduced the following generalization of $W_pH$: for some
convex, non-decreasing and left continuous function $\alpha$ on
$\rr^+$,

\beq \alpha(T_\VV(\nu,\mu)) \le 2C H(\nu|\mu), \ \ \forall \nu\in
M_1(\XX )\tag{$\alpha\!-\!T_\VV H(C)$} \neq and they established its
equivalence with some concentration inequality of the underlying
measure $\mu$ and of the i.i.d. sequences of common law $\mu$.

Recall that $T_\VV(\nu,\mu)=W_p(\nu,\mu)^2$ iff $\VV=\VV(p,d)$, the
family of all couples $(u,v)$ of real bounded measurable functions
on $\XX$ such that

\beq\label{A31} u(x) - v(y) \le d^p(x,y), \ \forall x,y\in E. \neq

Guillin-Leonard-Wu-Yao \cite{GLWY} propose a new
transport-information inequality, adapted to Markov processes (and in particular to consider deviation inequalities for integral functionals of Markov processes)

\beq W_p(\nu,\mu)^2 \le 4C^2 I(\nu|\mu), \ \ \forall \nu\in M_1(\XX
)\tag{$W_pI(C)$} \neq or the more general

\beq \alpha\left(T_\VV(\nu,\mu)\right) \le  I(\nu|\mu), \ \ \forall
\nu\in M_1(\XX )\tag{$\alpha\!-\! T_\VV I$}.\neq
 Using large deviations techniques they prove the
following characterization:

\bthm\label{thm-GLWY} {\rm (\cite{GLWY})} Let $((X_t)_{t\ge0},
\pp_\mu)$ be the $\mu$-symmetric and ergodic Markov associated with
the Dirichlet form $(\EE, \dd(\EE))$,  $\alpha: \rr^+\to
[0,+\infty]$ a left-continuous non-decreasing convex function with
$\alpha(0)=0$, and $\VV$ as above.

 The following properties are equivalent:

\bdes \item[(a)] $\mu$ satisfies the transport-information
inequality $(\alpha\!-\!T_\VV I)$.
\item[(b)] For all $(u,v)\in \VV$ and all $\lambda\ge 0$

\beq\label{L2bound} \lambda_{max}(\LL + \lambda u) :=
\sup_{g\in\dd(\EE): \mu(g^2)=1} \left[\lambda \int u g^2 d\mu -
\EE(g,g)\right] \le \lambda \mu(v)+ \alpha^*(\lambda) \neq where
$\LL$ is the generator of $(P_t)$ on $L^2(\XX ,\mu)$ and
$$\alpha^*(\lambda)=\sup_{r\ge0}\{\lambda r - \alpha (r)\}, \forall
\lambda\ge0
$$
is the semi-Legendre transformation of $\alpha$.

\item[(c)] For any initial measure $\nu=f\mu$ with $f\in L^2(\mu)$
and for  all $(u,v)\in \VV$

\beq\label{thmA11a} \pp_\nu\left(\frac 1t\int_0^t u(X_s)ds \ge
\mu(v)+ r\right)\le \|f\|_2 e^{- t\alpha(r)},\ \forall t,r>0. \neq

\ndes \nthm

\brmk{\rm The meaning of the deviation inequality characterization
(\ref{thmA11a}) of $\alpha\!-\!W_pI$ is clear in the ergodic
behavior of the Markov process $(X_t)$, as well as (\ref{L2bound})
in the study of the Schr\"odinger operator $\LL + u$. That is one
more reason why $\alpha\!-\!T_\VV I$ inequality is useful.}\nrmk

\brmk{\rm If $\VV$ is some family of $(u,u)\in (b\BB)^2$, (\ref{thmA11a})
becomes a deviation inequality of the empirical (time) mean from its
space mean $\mu(u)$ for the observable $u$ so that $(u,u)\in\VV$.
Notice that if $\VV=\{(u,u); u\in b\BB, \|u\|_{Lip}\le 1\}$ then
$T_\VV(\nu,\mu)=W_1(\nu,\mu)$, and $W_1I(C)$ is equivalent to the
Gaussian deviation inequality (\ref{thmA11a}) with
$\alpha(r)=r^2/(4C^2)$ for the Lipschitzian observable $u$ with
Lipschitzian coefficient $\|u\|_{Lip}\le 1$, which generalizes the
well known Hoeffding's inequality in the i.i.d. case. } \nrmk

Three criteria for $W_1I(C)$ are established in \cite{GLWY}:
spectral gap in $L^2(\mu)$; spectral gap in the space of Lipschitz
functions and a very general Lyapunov function criterion if
$\VV=\{(u,u); |u|\le \phi\}$ where $\phi>0$ is some fixed weight
funtion. And it is also shown that on a Riemannian manifold $\XX$
equipped with the Riemannian metric $d$, the log-Sobolev inequality

\beq H(\nu|\mu)\le 2C I(\nu|\mu),\ \forall \nu\in M_1(\XX ).
\tag{$HI(C)$} \neq implies $W_2I(C)$, which in turn implies the
Poincar\'e inequality

\beq Var_\mu(g)\le c_P \EE(g,g),\ \forall g\in L^2(\mu)\bigcap
\dd(\EE) \tag{$P(C)$} \neq where $Var_\mu(g)=\mu(g^2) -\mu(g)^2$ is
the variance. Furthermore $W_2I(C)\implies HI(C')$ once if the
Ricci-Bakry-Emery curvature of $\mu$ is bounded from below.

We organize this paper around the four questions below:

\bdes
\item[(i)] Investigate the relations between $W_pI$ with $W_pH$.
That is the objective of \S 2.

\item[(ii)] Prove that the spectral gap in the space of Lipschitz
functions implies a Cheeger type's isoperimetric inequality, which
is stronger than $W_1I$. That is the purpose of \S 3. We will also
establish deviation inequalities under natural quantities such as
the variance of the test function, refining \cite{GLWY}.

\item[(iii)] In \S4 we study relations between $(\alpha\!-\!W_2I)$ and the
$\beta$-log-Sobolev inequality:

\beq \label{Phi} \beta\circ\mu(g^2\log g^2)\le  \EE(g,g),\ \ \
\mu(g^2)=1, g\in\dd(\E),\neq where $\beta$ is a positive increasing
function. This inequality was connected in \cite{W07} to the well
developed $F$-Sobolev inequality introduced in \cite{W00}, so that
known criteria for the later can be applied directly to (\ref{Phi}).

\item[(iv)] Finally we present in \S 5 applications of $\Phi$-Sobolev inequality
$$
\|g^2\|_\Phi\le C_1 \EE(g,g) + C_2\mu(g^2)
$$
in transportation-information inequalities $\alpha\!-\!T_\VV I$ and
then in the concentration phenomena of $\frac 1t \int_0^t u(X_s)ds$
under integrability conditions on $u$.
 \ndes

\section{$W_pI$ implies $W_pH$ on a Riemannian manifold}

Recall (cf. Villani \cite{Vill03}) the well known Kantorovitch's
dual characterization :

\beq\label{A32} W_{p}^p(\nu, \mu)=\sup_{(u,v)\in\VV(p,d)} \int u
d\nu - \int v d\mu \neq where $\VV(p,d)$ is given in (\ref{A31}),
and Kantorovitch-Robinstein's identity

\beq\label{A33} W_{1}(\nu,\mu)=\sup_{\|u\|_{Lip}\le 1} \int u
d(\nu-\mu). \neq Throughout this section $\XX$ is a connected
complete Riemannian manifold equipped with the Riemannian metric
$d$, and $\mu=e^{-V} dx/Z$ ($Z$ being the normalization constant
assumed to be finite) with $V\in C^1(\XX )$, and $(\EE, \dd(\EE))$
is the closure of
$$
\EE(g,g):= \int_E |\nabla g|^2 d\mu(x), \ g\in C_0^\infty(\XX )
$$
where $\nabla$ is the gradient on $\XX$, and $C_0^\infty(\XX )$ is
the space of infinitely differentiable functions on $\XX$ with
compact support. In such case our Fisher information of $f\mu$ with
$0<f\in C^1(\XX )$ w.r.t. $\mu$ becomes
$$
I(f\mu|\mu)= \frac 14 \int \frac{|\nabla f|^2}{f} d\mu =  \frac 14
\int |\nabla \log f|^2 d\mu.
$$

\subsection{$W_1I(C)\implies W_1H(C)$}

\bthm\label{thm21} Assume that $\mu$ satisfies $W_1I(C)$. Then
$$
W_1(\nu,\mu)^2 \le 2C H(\nu|\mu),\ \forall \nu\in \MM_1(\XX )
$$
i.e., $\mu$ satisfies $W_1H(C)$.  \nthm

\bprf By Bobkov-G\"otze's criterion \cite{BG99} for $W_1H(C)$, it is
enough to show that for any bounded $g\in C^1(\XX )$ with $|\nabla
g|\le 1$ and $\lambda\ge 0$,

\beq\label{thm21b} \int e^{\lambda(g-\mu(g))}d\mu \le
e^{\lambda^2C^2/2}. \neq To this end we may assume that $\mu(g)=0$.
Consider

$$
Z(\lambda)=\int e^{\lambda g} d\mu, \ \mu_\lambda:= \frac{e^{\lambda
g}}{Z(\lambda)} \mu.
$$
We have by Kantorovitch's identity (\ref{A33})
$$
\frac{d}{d\lambda} \log Z(\lambda)= \mu_{\lambda} (g)\le
W_1(\mu_\lambda, \mu)
$$
but by $W_1I(C)$,

$$
W_1(\mu_\lambda, \mu)\le 2C \sqrt{I(\mu_\lambda|\mu)} = C\lambda
\sqrt{\int |\nabla g|^2 d\mu_\lambda}\le C\lambda.
$$
Thus
$$
\log Z(\lambda) \le \int_0^\lambda Ct dt = \frac{C\lambda^2}{2}
$$
the desired control (\ref{thm21b}).
 \nprf

The implication ``$W_1I(C)\implies W_1H(C)$" is strict, as shown by
the following simple counter-example (\cite{DGW03}).

\bexa {\rm Let $\XX=[-2, -1]\bigcup [1,2]$ and $\mu(dx)= (1_{[-2,
-1]} + 1_{[1,2]}) dx/2$. The Dirichlet form $(\EE, \dd(\EE))$ is
given by
$$
\EE(f,f)=\int f'^2 d\mu(x),\ \forall f\in\dd(\EE)=H^1(\XX )
$$
where $H^1(\XX )$ is the space of those functions $f\in L^2(\mu)$ so
that $f'\in L^2(\mu)$ (in the distribution sense). It corresponds to
the reflecting Brownian Motion in $\XX$, which is not ergodic. But
$W_1I(C)$ implies always the ergodicity. Thus $\mu$ does not satisfy
$W_1I(C)$. However $\mu$ satisfies $W_1H(C)$ by the Gaussian
integrability criterion in \cite{DGW03}.
 } \nexa

The argument above can be extended to more general transportation
information inequality $\alpha\!-\!W_1I$:

\bprop \label{prop21} Let $\alpha: \rr^+\to [0,+\infty]$ be a
left-continuous non-decreasing convex function with $\alpha(0)=0$.
Assume that $\mu$ satisfies $\alpha\!-\!W_1I$. Then $\mu$ satisfies

\beq\label{prop21a} \tilde \alpha(W_1(\nu,\mu))\le H(\nu|\mu),\
\forall \nu\in \MM_1(\XX ) \tag{$\tilde \alpha-W_1H$} \neq where
$\dsp \tilde\alpha(r)=2\int_0^r \sqrt{\alpha(s)} ds. $ In particular
for any Lipschitzian function $g$ with $\|g\|_{Lip}\le 1$,

$$
\mu(g> \mu(g) +r)\le e^{-\tilde \alpha(r)}, \ \forall r>0.
$$
 \nprop

\bprf By Gozlan-L\'eonard's criterion \cite{GL} for
$\tilde\alpha\!-\!W_1H$, it is enough to show that for any bounded
$g\in C^1(\XX )$ with $|\nabla g|\le 1$ and $\lambda\ge 0$,

\beq\label{prop21c} \int e^{\lambda(g-\mu(g))}d\mu \le e^{\tilde
\alpha^*(\lambda)}, \ \forall \lambda\ge0 \neq which implies the
last concentration inequality in this Proposition by Chebychev's
inequality. To show (\ref{prop21c}) we may assume that $\mu(g)=0$.
Let $Z(\lambda)$ and $\mu_\lambda$ be as in the previous proof of
Theorem \ref{thm21}, we have
$$
\frac{d}{d\lambda} \log Z(\lambda)= \mu_{\lambda} (g)\le
W_1(\mu_\lambda, \mu).
$$
But by the assumed $\alpha\!-\!W_1I$,

$$
W_1(\mu_\lambda, \mu)\le \alpha^{-1}\left( I(\mu_\lambda|\mu)\right)
= \alpha^{-1}\left(\frac{\lambda^2}4 \int |\nabla g|^2
d\mu_\lambda\right)\le \alpha^{-1}(\lambda^2/4)
$$
where $\alpha^{-1}(t):=\inf\{t\ge0;\ \alpha(r)>t\}$, $t\ge0$. Thus
$$
\log Z(\lambda) \le \int_0^\lambda \alpha^{-1}(t^2/4)
dt=:h(\lambda).
$$
Now by Fenchel-Legendre theorem, $h=(h^*)^*$, but
$$
h^*(r)=\sup_{\lambda\ge 0} (\lambda r - h(\lambda))=2\int_0^r
\sqrt{\alpha(s)} ds,
$$
which completes the proof of the desired control (\ref{prop21c}).
 \nprf

\subsection{$W_2I(C)\implies W_2H(C)$}

\bthm\label{thm22} Assume that $\mu$ satisfies $W_2I(C)$. Then
$$
W_2(\nu,\mu)^2 \le 2C H(\nu|\mu),\ \forall \nu\in \MM_1(\XX )
$$
i.e., $\mu$ satisfies $W_2H(C)$.  \nthm

\bprf We shall use the method of Hamilton-Jacobi equation due to
Bobkov-Gentil-Ledoux \cite{BGL01}. Consider the inf-convolution
$$
Q_t g(x):=\inf_{y\in E} (g(y)+ \frac 1{2t} d^2(x,y))
$$
which is viscosity solution of the Hamilton-Jacobi equation

\beq\label{HJequ}
\partial_t Q_t g + \frac 12 |\nabla Q_t g|^2 =0.
\neq By Bobkov-G\"otze's criterion \cite{BG99} for $W_2H(C)$, it is
enough to show that for any $g\in C^1_b(\XX )$,

\beq\label{thm22b} \int e^{Q_1g/C}d\mu \le e^{\mu(g)/C}. \neq To
this end we may and will assume that $\mu(g)=0$. Let
$\lambda=\lambda(t)=\kappa t$ where $\kappa>0$ will be determined
later and consider
$$
Z(t)=\int e^{\lambda Q_tg} d\mu, \ \mu_t:= \frac{e^{\lambda
Q_tg}}{Z(t)} \mu.
$$
We have

$$
\aligned \frac d{dt} \log Z(t)
&= \frac{1}{Z(t)}\int \left[
\lambda'(t) Q_t f+
\lambda(t) \partial_t Q_t g \right] e^{\lambda Q_tg} d\mu\\
&= \kappa \int Q_t g d\mu_t - \frac{\lambda}{2} \int |\nabla Q_t
g|^2
d\mu_t\\
&= \kappa \int Q_t g d\mu_t - \frac{2}{\lambda} I(\mu_t|\mu).
\endaligned
$$
But by Kantorovitch's identity (\ref{A32}),

\beq\label{thm22c} \int Q_t g d\mu_t \le \frac 1{2t} W_2^2 (\mu_t,
\mu) \neq and the assumed $W_2I(C)$ gives $W_2^2 (\mu_t, \mu)\le
4C^2 I(\mu_t|\mu)$. Thus for every $t>0$,
$$
\frac d{dt} \log Z(t)\le \left(\frac{2\kappa C^2}{t} -
\frac{2}{\kappa t} \right) I (\mu_t| \mu)
$$
Putting $\kappa=1/C$, we obtain $ \frac d{dt} \log Z(t)\le 0$ for
all $t>0$, which implies by the continuity of $\log Z(t)$ on $\rr^+$
that
$$
\int e^{Q_1g/C}d\mu = Z(1) \le Z(0)=1
$$
the desired (\ref{thm22b}).
 \nprf

\brmk {\rm The proof above is adapted from that of
Bobkov-Gentil-Ledoux \cite{BGL01} for the implication $HI(C)\implies
W_2H(C)$, originally established by Otto-Villani \cite{OVill00}. }
\nrmk

\brmk{\rm
We have thus established in this section
$$HI(C)\,\Rightarrow \,W_2I(C)\,\Rightarrow \, W_2H(C).$$
It was also established in \cite{GLWY} that under a lower bound of
the Ricci-Bakry-Emery curvature of $\mu$ that $W_2I$ implies back to
$HI$, and with additional conditions on this lower bound that $W_2H$
implies back $HI$. It is then a natural question to know if the
condition on the lower bound of the Ricci-Bakry-Emery curvature is
also necessary to get the reverse implication. A partial answer was
provided in \cite{p-CG04} where an example of a real probability
measure, with infinite lower bounded curvature, was shown to verify
$W_2H$ but not $HI$. Inspired by this example, we furnish here an
example where $W_1I$ holds (using Lyapunov conditions of
\cite[Section 5]{GLWY}) but not $HI$. Let then consider
$d\mu(x)=e^{-V(x)}dx$, where $V$ is symmetric $C^2$ (at least) and
given for large $x$ by
$$V(x)=x^4+4x^3\sin^2(x)+x^\beta.$$
Consider also the natural reversible process associated to this measure given by generator $Lf=f''-V'f'$.
Using $W(x)=e^{ax^4}$, by easy calculus, one sees that $LW\le -cx^4W+b$ (for some positive $b$ and $c$) if $\beta>2$. This Lyapunov condition also implies a Poincar\'e inequality (see \cite{BBCG} for example), so that using a slight modification of \cite[Lem. 5.7]{GLWY}, we get that $W_1I$ holds and also $W_2H$ by \cite{p-CG04}. Remark now that if $\beta<3$ then $V/V'^2$ is not bounded, which is a known necessary condition for $HI$ to hold (see \cite{p-CG04}). Unfortunatly, we are not up to now able to prove that $W_2I$ holds.
}\nrmk

\section{$W_1I$ and the isoperimetric inequality of Cheeger's type\\
 by means of the spectral gap in $C_{Lip}$}
In  this section we return to the general Polish space case $(\XX
,d)$. We assume that $\mu$ charges all non-empty open subsets of
$\XX$.

Let $C_{Lip}$ be the space of all real functions $g$ on $\XX$ which
are Lipschitz-continuous,  i.e., $\|g\|_{Lip}:=\sup_{x\ne y}\frac
{|f(x)-f(y)|}{d(x,y)}<+\infty$. We assume that there is an algebra
$\AA\subset C_{Lip}\bigcap \dd_2(\LL)$ (here $\dd_2(\LL)$ is the
domain of the generator $\LL$ in $L^2(\mu)$ associated with $(\EE,
\dd(\EE))$), which is a form core for $(\EE, \dd(\EE))$. Hence the
carr\'e-du-champs operator
$$
\Gamma(f,g):= \frac 12\left(\LL (fg) - f\LL g - g\LL f\right), \
\forall f,g\in \AA
$$
admits a unique continuous extension $\Gamma: \dd(\EE)\times
\dd(\EE)\to L^1(\XX , \mu)$. Throughout this section we assume that
$\Gamma$ is a differentiation, that is, for all $(h_k)_{1\le k\le
n}\subset \AA, g\in \AA$ and $F\in C_b^1(\rr^n)$,
$$
\Gamma(F(h_1,\cdots, h_n), g) = \sum_{i=1}^n
\partial_iF(h_1,\cdots, h_n)\Gamma(h_i, g)
$$
(this can be extended to $\dd(\EE)$).

\bthm\label{thm31} Assume that $\int d^2(x,x_0) d\mu(x)<+\infty$ for
some (or all) $x_0\in E$  and $\Gamma$ is a differentiation. Suppose
that there is a form core $\DD\subset C_{Lip}\bigcap \dd_2(\LL)$  of
$(\EE, \dd(\EE))$ such that $1\in\DD$ and

\beq\label{thm31a} W_1(\nu,\mu)=\sup_{g\in\DD: \|g\|_{Lip}\le 1}
\{\int g d(\nu-\mu)\}, \ \forall \nu \text{ with }
I(\nu|\mu)<+\infty \neq and

\beq\label{thm31b} \sqrt{\Gamma (g,g)} \le \sigma \|g\|_{Lip}, \
\mu-a.s., \ \forall g\in C_{Lip}\bigcap \dd_2(\LL)\neq and for some
constant $C>0$ and for any $g\in \DD$ with $\mu(g)=0$, there is
$G\in C_{Lip}\bigcap \dd_2(\LL)$ so that \beq\label{thm31c}  -\LL G
=g,\ \|G\|_{Lip}\le C \|g\|_{Lip}. \neq Then the Poincar\'e
inequality holds with $c_P\le C$, and the following isoperimetric
inequality of Cheeger's type

\beq\label{thm31d} W_1(f\mu, \mu)\le \sigma C \int
\sqrt{\Gamma(f,f)} d\mu,  0\le f\in \dd(\EE), \mu(f)=1\neq holds
true. In particular,

\beq\label{thm31e} W_1(\nu, \mu)^2\le 4(\sigma C)^2 I(\nu|\mu),\
\forall \nu\in\MM_1(\XX ). \neq Furthermore for any observable $g$
with $\|g\|_{Lip}=1$,

\beq\label{thm31f} \int g d(\nu-\mu) \le
2\sqrt{I(\nu|\mu)\left[\frac{V(g)}{2} + 2(\sigma C)^2 \sqrt{c_P
I(\nu|\mu)} \right]} \neq and for any $t,r,\delta>0$,

\beq\label{thm31g} \aligned &\pp_\beta\left(\frac 1t\int_0^t g(X_s)
ds
> \mu(g) + r\right)\\
&\le \|\frac{d\beta}{d\mu}\|_2 \exp\left(-t\frac{r^2}{(1+\delta)V(g)
+ \sqrt{[(1+\delta)V(g)]^2 + \frac{8c_P (\sigma C)^4}{\delta
V(g)}r^2 }}\right)\endaligned \neq where
$V(g):=\lim_{t\to\infty}\frac 1t Var_{\pp_\mu}\left(\int_0^t g(X_s)
ds\right)= 2 \int_0^\infty \<g-\mu(g), P_t g\>_\mu dt$ is the
asymptotic variance of $g$.
 \nthm

 \bprf Under the Lipschitzian spectral gap condition (\ref{thm31c}), it is noted in \cite{GLWY} that
 the Poincar\'e inequality holds with $c_P\le C$.

 For both (\ref{thm31d}) and (\ref{thm31e})  we may assume that $\nu=f\mu$ with $f\in \dd(\EE)$, $f\ge \vep>0$.
For any $g\in \DD$ with $\|g\|_{Lip}\le 1$ and $\mu(g)=0$, letting
$G:=(-\LL)^{-1}g$ be the unique solution of the Poisson equation
with $\mu(G)=0$, we have

$$
\int g d\nu - \int g d\mu= \<g,f\>_\mu= \EE(G, f)=\int \Gamma(G,f)
d\mu\le \|\sqrt{\Gamma(G,G)}\|_\infty \int \sqrt{\Gamma(f,f)} d\mu.
$$
Taking the supremum over all such $g$ and observing
$\|\sqrt{\Gamma(G,G)}\|_\infty\le \sigma \|G\|_{Lip}\le \sigma C$ we
obtain (\ref{thm31d}). Furthermore by Cauchy-Schwarz and the fact
that $\Gamma$ is a differentiation, we have
$$
\int \sqrt{\Gamma(f,f)} d\mu \le \sqrt{\int \frac{\Gamma (f,f)}{f}
d\mu} \sqrt{\int f d\mu} = 2\sqrt{I(\nu|\mu)}
$$
where (\ref{thm31e}) follows from (\ref{thm31d}).

For (\ref{thm31f}) writing $f=h^2$, we have
$$
\aligned \int g d\nu - \int g d\mu&=\int \Gamma(G,f) d\mu= 2\int h \Gamma(G,h) d\mu\\
&\le 2 \sqrt{\int \Gamma(h,h) d\mu\cdot\int \Gamma(G,G) h^2 d\mu}.
\endaligned$$
Using the inequality in \cite[Theorem 3.1]{GLWY}
$$
\int \Gamma(G,G) h^2 d\mu - \int \Gamma(G,G) d\mu \le
\|\Gamma(G,G)\|_\infty \|h^2\mu-\mu \|_{TV} \le (\sigma C)^2 \sqrt{4
c_P I(\nu|\mu)}
$$
and noting that $V(g)=2\<(-\LL)^{-1}g,g\>_\mu=2\EE(G,G)=2\int
\Gamma(G,G) d\mu$, we obtain

$$
\int g d\nu - \int g d\mu\le 2\sqrt{I(\nu|\mu)\left[\frac{V(g)}{2} +
2(\sigma C)^2 \sqrt{c_P I(\nu|\mu)} \right]}
$$
which is (\ref{thm31f}). Using $2 I^{3/2} \le \vep I + I^2/\vep$ in
(\ref{thm31f}), we obtain (\ref{thm31g}) by Theorem \ref{thm-GLWY}.
 \nprf

\brmk{\rm The $W_1I(\sigma C)$ inequality (\ref{thm31e}) is due to
Guillin and al. \cite{GLWY}, but the method therein is based on
the Lyons-Meyer-Zheng forward-backward martingale decomposition.
The argument here is simpler and direct, and yields the stronger
Cheeger type's isoperimetric inequality (\ref{thm31d}).
 } \nrmk

\brmk{\rm Letting $\delta$ be close to $0$, we see that
(\ref{thm31g}) is sharp for small $r$ by the central limit theorem.

Set $C_{Lip,0}=\{g\in C_{Lip,0};\ \mu(g)=0\}$. Under the
Lipschitzian spectral gap condition (\ref{thm31c}), the Poisson
operator $(-\LL)^{-1}: C_{Lip,0}\to C_{Lip,0}$ is a well defined
bounded linear operator w.r.t. the Lipschitzian norm, and the best
constant $C$ in (\ref{thm31c}) is the
 Lipschitzian norm $\|(-\LL)^{-1}\|_{Lip}$  and will
 be  denoted by $c_{Lip,P}$ (the index $P$ is referred to Poincar\'e).
}\nrmk

We now present four examples for illustrating usefulness of Theorem
\ref{thm31}.

 \bexa {\rm {\bf (Ornstein-Uhlenbeck process)} Consider the Ornstein-Uhlenbeck
 process $dX_t=\sqrt{2}dB_t -\sigma^{-2} X_t dt$ on $\XX=\rr$ where $\sigma>0$ and $B_t$ is the standard Brownian motion
 on $\rr$. Its
 unique invariant measure is $\mu=\NN(0,\sigma^2)$. For $f\in
 C_b^\infty(\rr)$, from the explicit solution $X_t= e^{-\sigma^{-2}t}\left(X_0+\int_0^te^{\sigma^{-2}s} \sqrt{2} dB_s
  \right)$, we see that $(P_tf)'= e^{-\sigma^{-2}t} P_tf'$. Hence
  $c_{Lip,P}=\|(-\LL)^{-1}\|_{Lip}=\sigma^2$. Therefore $\mu$ satisfies
  $W_1I(C)$ with $C=c_{Lip}=\sigma^2$ by Theorem \ref{thm31}.

  Furthermore $C=c_{Lip}=\sigma^2$ is also the best constant in
  $W_1I(C)$. Indeed $W_1I(C)\implies W_1H(C)$ and the best constant
  in $W_1H(C)$ of $\mu$ is $C=\sigma^2$.
  In other words Theorem \ref{thm31} produces the exact best
  constant $C$ in $W_1I(C)$ for this example.
}
 \nexa

\bexa{\bf (Reflected Brownian Motion)}  {\rm  Consider the reflected
Brownian Motion $X_t^D$ on the interval $\XX=[0, D]$ $(D>0)$
equipped with the usual Euclidean metric, whose generator is given
by $\LL f=f^{\prime\prime}$ with Neumann boundary condition  at
$0,D$. The unique invariant measure $\mu$ is the uniform law on $[0,
D]$. For every $g\in C_b^2([0, D])$ with $\int_{0}^{D} g(x)dx=0$,
the solution $G$ of the Poisson equation $-\LL G=g$ satisfies
$$
G'(x) = -\int_{0}^x g(t) dt, \ x\in [0, D].
$$
It is now easy to see that
$c_{Lip,P}=\sup_{\|g\|_{Lip}=1}\|G'\|_\infty$ is attained with
$g(x)=x-D/2$ and then $c_{Lip,P} =D^2/8$. Thus by Theorem
\ref{thm31}, the optimal constant $C_{W_1I}$ for this process
satisfies $C_{W_1I}\le c_{Lip}=D^2/8$. In comparison recall that the
best Poincar\'e constant $c_P=D^2/\pi^2$.

Since $W_1I(C)\implies W_1H(C)$ and the best constant of $W_1H(C)$
for the uniform law $\mu$ on $[0,D]$ is $D^2/12$, so we obtain
$$
\frac {D^2}{12} \le c_{W_1I} \le \frac{D^2}{8}.
$$
We do not know the exact value of $c_{W_1I}$ for this simple
example.
  } \nexa

\bexa{\rm Let $\XX$ be a compact connected Riemannian manifold of
dimension $n$ with empty or convex boundary. Assume that the Ricci
curvature is  nonnegative and its diameter is $D$. Consider the
Brownian Motion (with reflection in the presence of the boundary)
generated by the Laplacian operator $\Delta$.

 In \cite{Wu09} it is shown that $c_{Lip,P}=\|(-\Delta)^{-1}\|_{Lip}\le
 D^2/8$ (the latest quantity is exactly $c_{Lip,P}$ for the reflected Brownian Motion on the interval
 $[0,D]$). Thus by Theorem \ref{thm31}, $W_1I(C)$ holds with
 $C=D^2/8$.

 See \cite{Wu09} for more examples for which $c_{Lip,P}$ is
 estimated.
}\nexa

\bexa{\rm ({\bf One-dimensional diffusions}) Now let us consider the
one-dimensional diffusion with values in the interval $(x_0,y_0)$
generated by
$$
\LL f = a(x) f^{\prime\prime}+b(x) f',  f\in C_0^\infty(x_0,y_0)
$$
where $a,b$ are continuous such that $a(x)>0$ for all $x\in
(x_0,y_0)$. Let $((X_t)_{0\le t<\tau},\pp_x)$ be the martingale
solution associated with $\LL$ and initial position $x$, where
$\tau$ is the explosion time. With a fixec $c\in (x_0,y_0)$,
$$
s'(x):=\exp\left(-\int_c^x \frac {b(t)}{a(t)} dt\right),\ m'(x):=\frac 1{a(x)} \exp\left(\int_c^x \frac {b(t)}{a(t)} dt\right)
$$
 are respectively the derivatives of Feller's scale and speed
functions. Assume that
\begin{equation}\label{D1}
Z:=\int_{x_0}^{y_0} m'(x)\,dx<+\infty
\end{equation}
and let $\mu(dx)=m'(x)dx/Z$. It is well known that $(\LL,
C_0^\infty(x_0,y_0))$ is symmetric on $L^2(\mu)$.
\\
Assume also that
\begin{equation}\label{D2}
\int_c^{y_0} s'(x)\,dx \int_c^x m'(x)\,dx=\int_{x_0}^c
s'(x)\int_x^c m'(x)\,dx=+\infty
\end{equation}
which, in the Feller's classification, means that $x_0, y_0$ are
no accessible or equivalently $\tau=\infty$, $\pp_x$-a.s. In this
case by the $L^1$-uniqueness in \cite{Wu99}, the Dirichlet
form
\begin{eqnarray*}
\dd(\EE)&=&\left\{f\in\AA\CC(x_0, y_0)\bigcap L^2(\mu);
\int_{x_0}^{y_0}(f')^2d\mu<+\infty\right\}, \\
\EE(f,f)&=&\int_{x_0}^{y_0}(f')^2d\mu,\quad f\in\dd(\EE)
\end{eqnarray*}
is associated with $(X_t)$, where $\AA\CC(x_0, y_0)$ is the space
of absolutely continuous functions on $(x_0,y_0)$.

Fix some $\rho\in C^1(x_0,y_0)$ such that $\rho\in L^2(\mu)$ and
$\rho'(x)>0$ everywhere, consider the metric
$d_\rho(x,y)=|\rho(x)-\rho(y)|$. A function $f$ on $(x_0,y_0)$ is
Lipschitz with respect to $d_\rho$ (written as $f\in
C_{\textrm{Lip}(\rho)}$) if and only if $f\in\AA\CC(x_0,y_0)$ and
$$
\|f\|_{\textrm{Lip}(\rho)}=\sup_{x_0<x<y<y_0}
\frac{|f(y)-f(x)|}{\rho(y)-\rho(x)}=\|\frac {f'}{\rho'}\|_\infty.
$$
The argument below is borrowed from \cite{p-DMW}. Assume that

\begin{equation}\label{C-rho}
C(\rho):=\sup_{x\in (x_0,y_0)} \frac 1{\rho'(x)} \int_x^{y_0}
[\rho(t)-\mu(\rho)] m'(t)\,dt<+\infty.
\end{equation}
For every $g\in C_{\textrm{Lip}(\rho)}$ with $\mu(g)=0$, then $
f(x)=\int_c^x dy \int_y^{y_0} g(t) m'(t)\,dt - A
$ (in $C^2$) solves
\begin{equation}\label{cor44a}
-(a f^{\prime\prime} + b f') = g.
\end{equation}
It is obvious that
$$
\|f\|_{\textrm{Lip}(\rho)}=\sup_{x\in (x_0,y_0)} \frac 1{\rho'(x)}
\int_x^{y_0} g(t) m'(t)\,dt.
$$
An elementary exercise (as done in \cite{p-DMW}) shows that the
last quantity is always not greater than
$C(\rho)\|g\|_{\textrm{Lip}(\rho)}$. Thus $f\in L^2(\mu)$ (for
$\rho\in L^2(\mu)$). By Ito's formula, $f\in \dd_2(\LL)$.   With
the constant $A$ so that $\mu(f)=0$, $f$ given above is the unique
solution in $L^2(\mu)$ with zero mean of (\ref{cor44a}) by the
ergodicity of $(X_t)$.
 We see also
that $C(\rho)$ is the best constant by taking $g=\rho-\mu(\rho)$. In
other words condition (\ref{thm31c}) is verified with the best
constant $C=c_{Lip,P}=C(\rho)$. Hence from Theorem \ref{thm31}, we
get

\bcor\label{cor44} Let $a,b: (x_0,y_0)\to \rr$ be continuous such
that $a(x)>0$ for all $x$ and conditions (\ref{D1}) (\ref{D2}) be
satisfied. Assume (\ref{C-rho}) and $\sigma:=\sup_{x\in (x_0,y_0)}
\sqrt{a(x)} \rho'(x)<+\infty$. Then $\mu$ satisfies $W_1I(\sigma
C(\rho))$ on $((x_0,y_0), d_{\rho})$. In particular for
$$
\rho_a(x)=\int_c^x \frac 1{\sqrt{a(t)}} dt
$$
($d_{\rho_a}$ is the metric associated with the carr\'e-du-champs
operator of the diffusion), if $C(\rho_a)<+\infty$, then $\mu$
satisfies $W_1I(C(\rho_a))$ on $((x_0,y_0), d_{\rho_a})$. \ncor

\brmk{\rm The quantity $C(\rho)$ in (\ref{C-rho}) is not innocent:
Chen-Wang's variational formula for the spectral gap $\lambda_1$
says that (\cite{Chen05, Wang}): $\dsp \lambda_1=\sup_{\rho} \frac
1{C(\rho)}. $ } \nrmk

}\nexa

\section{Functional inequalities and $W_2I$ inequalities}

Throughout this section we consider the framework of Section 2, i.e.
$\XX$ is a connected complete Riemannian manifold $M$ with $\mu(\d
x):= \e^{-V(x)}\d x/Z$ for some $V\in C(M)$ with $Z:= \int_M
\e^{-V(x)}\d x<\infty.$ Recall that in \cite{OVill00, GLWY} was proven the fact that a logarithmic Sobolev inequality implies $W_2I$, and that (using HWI inequalities) under a lower bounded curvature, the converse was also true. We extend here this assertion for $\alpha-W_2I$ inequalities.

\beg{thm} \label{T4.1}  $(1)$ Let $\beta\in C([0,\infty))$ be
increasing with $\beta(0)=0$ such that

$$\gamma(r):= \ff 1 2 \int_0^r \ff{\d s}{\ss{\beta(s)}}<\infty,\ \ \ \
r>0.$$ Then the following $\beta$-log-Sobolev inequality

\beq\label{P1} \beta\circ \mu(g^2\log g^2)\le \mu(|\nnn g|^2),\ \ \
g\in C_b^1(M), \mu(g^2)=1,\end{equation} implies

\beq\label{AA} \aa(W_2(\nu,\mu))\le I(\nu|\mu),\ \ \ \ \nu\in
M_1(\XX )\neq for $\aa(s):= \beta\circ \gamma^{-1}(s), s\ge 0.$

$(2)$ Assume that $\Ric +\Hess_V\ge -K$ for some $K\ge 0.$ Then
$(\ref{AA})$ implies $(\ref{P1})$ for

$$\beta(r):= \inf\big\{s>0:\ 2\ss{2s}\aa^{-1}(s) +K (\aa^{-1}(s))^2\ge
r\big\},\ \ \ r\ge 0.$$
\end{thm}

\beg{proof} (1) According to \cite[Theorem 2.2]{W07}, (\ref{P1})
implies

\beq\label{TP} W_2(f\mu, \mu)\le \gamma\circ\mu(f\log f),\ \ \ f\ge
0, \mu(f)=1.\end{equation} Then (\ref{AA}) follows from (\ref{TP})
and (\ref{P1}). For readers' convenience, we include below a brief
proof of (\ref{TP}), inspired by the seminal work \cite{OVill00}  pushed further in \cite{W04}.

Since a continuous function can be uniformly approximated by
smooth ones,   we may and do assume that $V$ is smooth. Let $P_t$
be the diffusion semigroup generated by $\LL:= \Delta-\nnn V.\nnn$.
Then $P_t$ is symmetric in $L^2(\mu).$ For fixed $f>0$ with
$\mu(f)=1,$ let $\mu_t= (P_t f)\mu,\ t>0.$ According to \cite[page
176]{W04} for $p=2$ (see also \cite{OVill00} under a curvature
condition), we have

\beq\label{W2.1} \ff{\d^+}{\d t}\big\{-W_2(\mu,\mu_t)\big\} :=
\limsup_{s\downarrow 0} \ff{W_2(\mu,\mu_t)-W_2(\mu,\mu_{t+s})}s \le
2 \mu\big(\big|\nnn\ss{P_t f}\big|^2\big)^{1/2}.\end{equation} Let

$$\gamma(r)=\ff 1 2 \int_0^r\ff{\d
s}{\ss{\beta (s)}},\ \ \ r>0.$$ It suffices to prove for the case
that $\gamma(r)<\infty$ for $r>0.$ By (\ref{P1}) we have

\beg{equation*}\beg{split} &\ff{\d}{\d t} \gamma\circ \mu(P_t f\log
P_t f)= 4 \gamma'\circ \mu(P_t f\log P_t f) \mu\big(\big|\nnn
\ss{P_t
f}\big|^2\big)\\
&= \ff{2\mu\big(\big|\nnn \ss{P_t f}\big|^2\big)}{\ss{\beta
\circ\mu(P_t f\log P_t f)}} \ge 2 \mu\big(\big|\nnn \ss{P_t
f}\big|^2\big)^{1/2}.\end{split}\end{equation*} Combining this with
(\ref{W2.1}) we obtain

$$\ff{\d^+}{\d t}\big\{-W_2(\mu,\mu_t)\big\}\le \ff{\d}{\d t}
\gamma\circ \mu(P_t f\log P_t f),$$ which implies (\ref{TP}) by
noting that $P_t f\to \mu(f)=1$ as $t\to\infty$.

(2) By the HWI inequality (see \cite{OVill00,BGL01}), we have

$$\mu(g^2\log g^2) \le \ff{2(\e^{2Kt}-1)} K \mu(|\nnn g|^2)
+\ff{K\e^{2Kt}}{\e^{2Kt}-1} W_2(g^2\mu,\mu)^2,\ \ \ \mu(g^2)=1,
t>0.$$ Combining this with (\ref{AA}) we obtain

$$\mu(g^2\log g^2) \le \inf_{t>0}\Big\{\ff{2(\e^{2Kt}-1)} K \mu(|\nnn
g|^2) + \ff{K\e^{2Kt}}{\e^{2Kt}-1} [\aa^{-1} \big(\mu(|\nnn
g|^2)\big)]^2\Big\}.$$ Taking $t>0$ such that

$$\e^{2Kt}= 1 +K \ff{\aa^{-1}(\mu(|\nnn g|^2))}{ \ss{2\mu(|\nnn
g|^2)}},$$ we obtain

$$\mu(g^2\log g^2)\le 2 \ss{2 \mu(|\nnn g|^2)}\aa^{-1}(\mu(|\nnn
g|^2)) +K [\aa^{-1}(\mu(|\nnn g|^2))]^2.$$ This completes the proof.
\end{proof}

Let us give a natural family of examples, namely when $\beta$ is a power function.

\beg{cor}\label{C4.2} For any $\delta \in [1,2)$,

$$\mu(g^2\log g^2)^\delta \le C\mu(|\nnn g|^2),\ \ \ g\in C_b^1(\XX ),
\mu(g^2)=1$$ implies

$$W_2(\nu,\mu)^2 \le
\ff{C^{2/\delta}}{(2-\delta)^{2}}I(\nu |\mu)^{(2-\delta)/\delta}.$$
Inversely if $\Ric +\Hess_V$ is bounded below, then

$$W_2(\nu,\mu)^2 \le CI(\nu|\mu)^{(2-\delta)/\delta}$$ implies

$$\mu(g^2\log g^2)^\delta\le C' \mu(|\nnn g|^2),\ \ \ g\in
C_b^1(\XX ),\mu(g^2)=1$$ for some $C'>0.$ \end{cor}

\beg{proof} For $\beta(r):= r^\delta/C$ we have $\gamma(r)= \ff{\ss
C}{2-\delta} r^{(2-\delta)/2}$ so that

$$\beta\circ \gamma^{-1}(s)=\ff 1 C \Big(\ff{2-\delta}{\ss
C}s\Big)^{2\delta/(2-\delta)}=
\ff{(2-\delta)^{2\delta/(2-\delta)}}{C^{2/(2-\delta)}}s^{2\delta/(2-\delta)}.$$
Then the first assertion follows from Theorem \ref{T4.1}(1).

Next, for $\aa(r)= r^{2\delta/(2-\delta)} C^{-\delta/(2-\delta)},$
we have $\aa^{-1}(s)=\sqrt{C}s^{(2-\delta)/2\delta}.$ Since
$2-\delta\le 1,$ Theorem \ref{T4.1}(2) implies

\beq\label{W*}   \mu(g^2\log g^2)\le 2\ss{2C} \mu(|\nnn
g|^2)^{1/\delta} +KC \mu(|\nnn g|^2)^{(2-\delta)/\delta},\ \ \
\mu(g^2)=1.\end{equation}
 Since $2-\delta\le 1$, this implies

\beq\label{W*2}  \mu(g^2\log g^2)\le C' \mu(|\nnn g|^2)^{1/\delta},\
\ \ \mu(g^2)=1, \mu(|\nnn g|^2)\ge 1\end{equation} for some $C'>0.$
Moreover, since $\delta\ge 1,$ (\ref{W*}) implies the defective
log-Sobolev inequality

$$\mu(g^2\log g^2)\le C_1\mu(|\nnn g|^2)+C_2,\ \ \ \mu(g^2)=1$$ for
some $C_1, C_2>0,$ which in particular implies that the spectrum of
$L:= \Delta +\nnn V$ is discrete (see e.g. \cite{W00, Wu00b}), and
hence the Poincar\'e inequality holds since $\lambda_0=0$ is the
simple eigenvalue due to the connection of the manifold. Thus, the
strict log-Sobolev inequality

$$\mu(g^2\log g^2)\le C' \mu(|\nnn g|^2),\
\ \ \mu(g^2)=1$$ for some constant $C'>0.$ The proof is then
completed by combining this with (\ref{W*2}). \end{proof}

\bexa {\rm Let $\Ric$ be bounded below, and $\rho_o$ the Riemannian
distance function to a fixed point $o\in E.$ Let $V\in C(\XX )$ such
that $V-a \rho_o^\theta$ is bounded for some $a>0$ and $\theta\ge
2.$ Then (\ref{AA}) holds for $\aa(r)=Cr^{2(\theta-1)}$ for some
$C>0$, i.e.

\beq\label{W*3} C W_2(\nu,\mu)^{2(\theta -1)}\le I(\nu |\mu),\ \ \
\nu\in M_1(\XX ).\end{equation} The power $2(\theta-1)$ is sharp,
i.e. the above inequality does not hold if this power is replaced by
any larger number, as seen from Proposition \ref{prop21}.

Indeed, by \cite[Corollaries 2.5 and 3.3]{W00}, we have

$$\mu(g^2\log^{2(\theta-1)/\theta}(g^2 +1)) \le C_1\mu(|\nnn g|^2)
+C_2,\ \ \ \mu(g^2)=1$$ holds for some $C_1, C_2>0.$ By Jensen's
inequality we obtain

$$\mu(g^2\log g^2)^{2(\theta -1)/\theta}\le \mu(g^2\log
(g^2+1))^{2(\theta -1)/\theta} \le C_1\mu(|\nnn g|^2) +C_2,\ \ \
\mu(g^2)=1.$$ Combining this with the log-Sobolev inequality as in
the proof of Corollary \ref{C4.2}, we obtain

 $$\mu(g^2\log g^2)^{2(\theta -1)/\theta}\le \mu(g^2\log
(g^2+1))^{2(\theta -1)/\theta} \le C'\mu(|\nnn g|^2),\ \ \
\mu(g^2)=1$$ for some constant $C'>0.$ According to Corollary
\ref{C4.2}, this implies (\ref{W*3}).

} \nexa






\section{$\Phi$-Sobolev inequality and concentration inequality \\
for unbounded observables under integrability condition}

Let $\Phi:\rr^+\to [0,+\infty]$ be a convex, increasing and left
continuous function with $\Phi(0)=0$, such that

\beq\label{51} \lim_{r\to +\infty} \frac {\Phi(r)}{r}=+\infty. \neq
Consider the Orlicz space $L^{\Phi}(\mu)$ of those measurable
functions $g$ on $\XX$ so that its gauge norm

$$
N_{\Phi}(g):=\inf\{c>0; \int \Phi(|g|/c) d\mu\le 1\}
$$
is finite, where the convention $\inf \emptyset:=+\infty$ is used.
The Orlicz norm of $g$ is defined by
$$
\|g\|_\Phi:=\sup\{\int g u \,d\mu; \ N_{\Psi}(u)\le 1\}
$$
where

\beq \Psi(r):=\sup_{\lambda\ge 0} (\lambda r - \Phi(\lambda)),\
r\ge0 \neq is the convex conjugation of $\Phi$.  The so called
(defective) $\Phi$-Sobolev inequality says that for some two
nonnegative constants $C_1, C_2\ge 0$

\beq\label{52} \|g^2\|_{\Phi} \le C_1 \EE(g,g) + C_2\mu(g^2), \
\forall g\in \dd(\EE), \mu(g^2)=1. \neq Under the assumption of the
Poincar\'e inequality with the best constant $C_P$, (\ref{52}) can
be transformed into the following tight version

\beq \label{53} \|(g-\mu(g))^2\|_{\Phi} \le (C_1 + C_2C_P) \EE(g,g),
\ \forall g\in\dd(\EE) \neq called sometimes Orlicz-Poincar\'e
inequality.

\bthm\label{thm51} Assume the $\Phi$-Sobolev inequality (\ref{52})
and the Poincar\'e inequality with constant $C_P$. Then

\bdes \item[(a)] for any $\mu$-probability density $f$,

\beq \label{thm51a} \|f-1\|_\Phi \le \sqrt{C_1' I(f\mu|\mu)^2 +
C_2'I(f\mu|\mu)}\neq where $C_1'=(C_1+2C_2C_P)C_1,
C_2'=(C_1+2C_2C_P) \cdot 4C_2$; or equivalently for any observable
$u\in L^\Psi(\mu)$ ($\Psi$ being the convex conjugation of $\Phi$
given above) so that $N_\Psi(u)\le 1$ and for all $t,r>0$,

\beq\label{thm51b} \pp_\nu\left(\frac 1t\int_0^t u(X_s) ds >\mu(u) +
r\right)\le \|\frac {d\nu}{d\mu}\|_2\exp\left(- t \cdot\frac
{\sqrt{4C_1'r^2 + (C_2')^2}-C_2'}{2C_1'}\right). \neq

\item[(b)] for any $\mu$-probability density $f$,

\beq \label{thm51c}\sup_{u\in b\BB: N_\Psi(u^2)\le 1} \int (f-1) u
d\mu \le \sqrt{ 2 (C_1+4C_2C_P)I(f\mu|\mu)} \neq or equivalently for
any $u\in L^1(\mu)$ such that $u^2\in L^\Psi(\mu)$,

\beq\label{thm51d} \pp_\nu\left(\frac 1t\int_0^t u(X_s) ds >\mu(u) +
r\right)\le \|\frac {d\nu}{d\mu}\|_2\exp\left(- t \frac{r^2}{2
(C_1+4C_2C_P)\|u^2\|_\Psi} \right), \ \forall t,r>0. \neq

\item[(c)] More generally for any $p\in [1,+\infty)$, there is a constant
$\kappa>0$ depending only of $p, C_1, C_2, C_P$ such that for any
$\mu$-probability density $f$,

\beq \label{thm51e}\alpha_p\left(\sup_{u\in b\BB: N_\Psi(u^2)\le 1}
\int (f-1) u d\mu\right) \le I(f\mu|\mu) \neq where
$\alpha_p(r)=(1+r^2/\kappa )^{p/2}-1$; or equivalently for any $u\in
L^1(\mu)$ such that $N_\Psi(|u|^p)\le 1$,

\beq\label{thm51f} \pp_\nu\left(\frac 1t\int_0^t u(X_s) ds >\mu(u) +
r\right)\le \|\frac {d\nu}{d\mu}\|_2\exp\left(- t [1+r^2/\kappa
)^{p/2}-1] \right), \ \forall t,r>0. \neq
 \ndes \nthm

 As there are numerous practical criteria for the $\Phi$-Sobolev
 inequality (see e.g. \cite{Chen05, Led01, Wang}), this theorem is very useful and gives different
 concentration behaviors for $\frac 1t\int_0^t u(X_s) ds$,  according
 to the integrability condition $|u|^p\in L^\Psi(\mu)$ where $p\in
 [1,+\infty)$.

 This result generalizes the sharp concentration inequality under the log-Sobolev inequality
 in Wu \cite{Wu00c}. For applications of $\Phi$-Sobolev inequalities
 in large deviations see Wu and Yao \cite{WY07}.

\brmk{\rm  As the l.h.s. of (\ref{thm51a}), (\ref{thm51c}) and
(\ref{thm51e}) are the transportation cost $T_\VV (f\mu, \mu)$, with
$\VV=\{(u,u);\ u\in b\BB, \ N_\Psi(|u|^p)\le 1\}$, $p=1,2, p\ge 1$
respectively, so they are the transportation-information inequality.
In this point of view, the equivalence between (\ref{thm51a}) and
(\ref{thm51b}) in part (a), that between (\ref{thm51c}) and
(\ref{thm51d}) in part (b) and that between (\ref{thm51e}) and
(\ref{thm51f}) in part (c) are all immediate from Theorem
\ref{thm-GLWY} (the passage from bounded $u$ to general $u$ in the
concentration inequalities (\ref{thm51b}), (\ref{thm51d}) and
(\ref{thm51f}) can be realized easily by dominated convergence). }
\nrmk

\brmk {\rm  The concentration inequalities (\ref{thm51b}),
(\ref{thm51d}) and (\ref{thm51f}) are all sharp in order. Indeed
consider the Ornstein-Uhlenbeck process on $\rr$ generated by $\LL
f=f^{\prime\prime} - x f'$: the $\Phi$-Sobolev inequality $(X_t)$
holds with $\Phi(r)=(1+r)\log (1+r)$ and $\mu=\NN(0,1)$. Consider
$u(x):=|x|^{2/p}$ where $p\ge 1$. Then $u^p\in L^\Psi(\mu)$, and
$\frac 1t\int_0^t u(X_s) ds= \frac 1 t\int_0^t |X_s|^p ds$  possess
exactly the concentration behaviors exhibited by the r.h.s. of
(\ref{thm51f}) for large deviation value $r$, and for small
deviation value $r$ of order $1/\sqrt{t}$ if $t$ is large enough (by
the central limit theorem).
 }
 \nrmk

\bprf[Proof of Theorem \ref{thm51}] As explained in the previous
remarks, (\ref{thm51b}) (resp. (\ref{thm51d}); (\ref{thm51f})) is
equivalent to (\ref{thm51a})(resp.  (\ref{thm51c}), (\ref{thm51e})),
all by Theorem \ref{thm-GLWY}.\\
It is not surprising that the proof relies on the ideas first used in \cite{BV}, establishing criterions for $W_1H$ under integrability criteria. Note also that the reader may easily adapt the proof to use conditions on $F$-Sobolev inequalities (equivalent to some Orlicz-Poincar\'e inequality) and integrability on $u$ (rather than Orlicz norm of $u$).

(a) For  (\ref{thm51a}) we may assume that $I(f\mu|\mu)$ is finite,
i.e., $\sqrt{f}\in\dd(\EE)$ (and then $I(f\mu|\mu)=\EE(\sqrt{f},
\sqrt{f})$). For any $u\in L^\Psi(\mu)$ with $N_\Psi(u)\le 1$, we
have by Cauchy-Schwartz
$$
\aligned \int |(f-1) u| d\mu &= \int |\sqrt{f} -1|  (\sqrt{f} +1)
|u|
d\mu\\
&\le \sqrt{\int (\sqrt{f} -1)^2 |u| d\mu} \sqrt{\int  (\sqrt{f}
+1)^2
|u| d\mu }\\
&\le \sqrt{\|(\sqrt{f} -1)^2\|_\Phi \|(\sqrt{f} +1)^2\|_\Phi}
\endaligned
$$
But by the assumed $\Phi$-Sobolev inequality (\ref{52}),
$$
\|(\sqrt{f} -1)^2\|_\Phi\le C_1\EE(\sqrt{f}, \sqrt{f}) +C_2 \int
(\sqrt{f} -1)^2 d\mu
$$
and $\int (\sqrt{f} -1)^2 d\mu= 2(1-\mu(\sqrt{f}))\le 2
Var_\mu(\sqrt{f})\le 2 C_P \EE(\sqrt{f}, \sqrt{f})$; moreover
$$
\|(\sqrt{f} +1)^2\|_\Phi\le C_1 \EE(\sqrt{f}, \sqrt{f}) +C_2 \int
(\sqrt{f} +1)^2 d\mu
$$
and $\int (\sqrt{f} +1)^2 d\mu\le 4$. Thus we get
$$
\int |(f-1) u| d\mu \le \sqrt{(C_1+2C_2C_P)\EE(\sqrt{f}, \sqrt{f})
(C_1 \EE(\sqrt{f}, \sqrt{f}) + 4 C_2 )}
$$
where (\ref{thm51a}) follows by recalling $I(f\mu|\mu)=\EE(\sqrt{f},
\sqrt{f})$.

\vskip15pt (b) For any $u$ so that $N_\Psi(u^2)\le 1$ we use now
differently Cauchy-Schwartz inequality to get:
$$
\aligned \int |(f-1) u| d\mu &\le \sqrt{\int (\sqrt{f} -1)^2  d\mu}
\sqrt{\int  (\sqrt{f} +1)^2
u^2 d\mu }\\
\endaligned
$$
But as noticed in the proof of (a),
$$
\int (\sqrt{f} -1)^2  d\mu\le 2 Var_\mu(\sqrt{f})\le 2\min\{C_P
I(f\mu|\mu), 1\}
$$
and
$$
\int  (\sqrt{f} +1)^2 u^2 d\mu\le \|(\sqrt{f} +1)^2\|_\Phi\le C_1
I(f\mu|\mu)+ 4C_2.
$$
Plugging those two estimates into the previous inequality we get
(\ref{thm51c}).

\vskip15pt (c). Letting $q:=p/(p-1)$ we have by H\"older's
inequality,
$$
\aligned
 \int |(f-1) u| d\mu &\le (\mu(|f-1|))^{1/q}
\left(\int  |f-1| |u|^p d\mu \right)^{1/p}\endaligned$$ Note that
$\mu(|f-1|)\le 2$ and by \cite[Theorem 3.3]{GLWY},

$$(\mu(|f-1|))^2\le 4Var_\mu(\sqrt{f})\le 4C_P I$$ where $I:=I(f\mu|\mu).$
On the other hand by part (a),
$$
\int  |f-1| |u|^p d\mu\le \|f-1\|_\Phi \le  \sqrt{C_1' I^2 + C_2'I}.
$$
Substituting those estimates into the first inequality we get

$$
\aligned \left( \int |(f-1) u| d\mu\right)^2 &\le (\max\{4,
4C_PI\})^{1/q} \left( C_1' I^2 + C_2'I \right)^{1/p}\\
&\le \begin{cases} 4^{1/q}(C_1'+C_2'C_P)^{1/p} \cdot I^{2/p}, \ &\text{if } C_PI\ge 1; \\
4^{1/q}(C_1'+C_2'C_P)^{1/p} C_P^{(p-2)/p} \cdot I, &\text{
otherwise.}
\end{cases}\endaligned$$
The last term is less than $\kappa [(1+I)^{2/p}-1]$ for some
constant $\kappa>0$. That yields to (\ref{thm51e}).
 \nprf

Let us finally relate previous inequalities to usual $\alpha-WI$ inequalities.

\bcor\label{cor51} Assume the $\Phi$-Sobolev inequality (\ref{52})
and the Poincar\'e  inequality. Assume that $d^p(x,x_0)\in
L^\Psi(\mu)$ for some $p\ge1$ where $\Psi$ is the convex conjugation
of $\Phi$. Then there are positive constants $C_1', C_2'$ and
$\kappa$ such that for all $\nu\in \MM_1(\XX)$,
$$
W_p^p(\nu,\mu) \le \sqrt{C_1'I(\nu,\mu)^2 + C_2' I(\nu|\mu)},\
$$
and
$$
\kappa\left([1+W_1(\nu,\mu)^2]^{p/2}-1\right)\le  I(\nu|\mu)
$$
and when $p\ge 2$,
$$
\kappa\left([1+W_2(\nu,\mu)^4]^{p/4}-1\right)\le  I(\nu|\mu).
$$
 \ncor
\bprf Recall the following fact (\cite[Proposition 7.10]{Vill03}),
$$
W_p^p(\nu,\mu)\le 2^{p-1}\|d(\cdot, x_0)^p(\nu-\mu)\|_{TV}.
$$
Then this corollary follows directly from Theorem \ref{thm51}. \nprf


\end{document}